\magnification=\magstep1
\input amstex
\documentstyle{amsppt}

\define\defeq{\overset{\text{def}}\to=}
\define\ab{\operatorname{ab}}

\define\Gal{\operatorname{Gal}}
\define\diag{\operatorname{diag}}

\def \isom {\overset \sim \to \rightarrow}

\define\Pic{\operatorname{Pic}}
\define\Br{\operatorname{Br}}

\define\period{\operatorname{period}}
\define\index{\operatorname{index}}

\define\id{\operatorname{id}}
\define\cn{\operatorname{cn}}

\define\Ker{\operatorname{Ker}}

\def \c{\operatorname {c}}

\def \res{\operatorname {res}}
\def \Jac{\operatorname {Jac}}
\def \Ind{\operatorname {Ind}}

\def \et{\operatorname {et}}

\def\Im{\operatorname{Im}}

\def\and{\operatorname{and}}

\def\Res{\operatorname{Res}}
\def\Div{\operatorname{Div}}
\def\Coker{\operatorname{Coker}}
\def\cl{\operatorname{cl}}
\def\tor{\operatorname{tor}}

\NoRunningHeads
\NoBlackBoxes
\topmatter

\title
Arithmetic of $p$-adic curves and sections of geometrically abelian fundamental groups
\endtitle

\author
Mohamed Sa\"\i di
\endauthor

\abstract
Let $X$ be a proper, smooth, and geometrically connected curve of genus $g(X)\ge 1$ over a $p$-adic local field. We prove that there exists an effectively computable open affine
subscheme $U\subset X$ with the property that $\period (X)=1$, and $\index (X)$ equals $1$ or $2$ (resp. $\period(X)=\index (X)=1$, assuming $\period (X)=\index (X)$),  if 
(resp. if and only if) the exact sequence of the geometrically abelian fundamental 
group of $U$ {\it splits}. We compute the torsor of splittings of the exact sequence of the geometrically abelian absolute Galois group associated to $X$, and give a
new characterisation of sections of arithmetic fundamental groups of curves over $p$-adic local fields which are orthogonal to $\Pic^0$ (resp. $\Pic^{\wedge}$). As a consequence we observe that the {\it non-geometric} (geometrically pro-$p$) section constructed by Hoshi in [Hoshi] is orthogonal to $\Pic^0$.
\endabstract

\endtopmatter

\document

\subhead
\S 0. Introduction/Main Results
\endsubhead
Let $k$ be a field of characteristic $0$ and $X$ a proper, smooth, and geometrically connected curve
over $k$ of genus $g(X)\ge 1$ with function field $K\defeq k(X)$. Let $\eta$ be a geometric point of $X$ with values in its generic point. Thus, $\eta$
determines an algebraic closure $\overline K$ (resp. $\overline k$) of $K$ (resp. $k$). Let $U\subseteq X$ be a
non-empty open subscheme and $U_{\overline k}\defeq U\times_k\overline k$. We have an exact sequence of fundamental groups $1\to \pi_1(U_{\overline k},\eta)\to \pi_1(U,\eta)
\to G_k\defeq \Gal (\overline k/k)\to 1$ (here $\eta$ is the geometric point of $U$, $U_{\overline k}$, naturally induced by $\eta$). By pushing this sequence by the maximal abelian quotient 
$\pi_1(U_{\overline k},\eta)\twoheadrightarrow  \pi_1(U_{\overline k},\eta)^{\ab}$ of $\pi_1(U_{\overline k},\eta)$
we obtain an exact sequence 
$$1\to \pi_1(U_{\overline k},\eta)^{\ab}\to \pi_1(U,\eta)^{(\ab)}\to G_k\to 1,\tag 1$$ 
where $\pi_1(U,\eta)^{(\ab)}\defeq \pi_1(U,\eta)/\Ker (\pi_1(U_{\overline k},\eta)\twoheadrightarrow  \pi_1(U_{\overline k},\eta)^{\ab})$ is the {\it geometrically abelian} fundamental group of $U$. 
Similarly, by pushing the exact sequence of absolute Galois groups 
$1\to G_{\overline k(X)}\defeq 
\Gal (\overline K/\overline kK)\to G_{k(X)}\defeq \Gal (\overline K/K)\to G_k\to 1$ 
by the maximal abelian quotient
$G_{\overline k(X)}\twoheadrightarrow G_{\overline k (X)}^{\ab}$ of $G_{\overline k(X)}$ we obtain an exact sequence
$$1\to G_{\overline k(X)}^{\ab}\to G_{k(X)}^{(\ab)}\to G_k\to 1,\tag 2$$ 
where $G_{k(X)}^{(\ab)}\defeq G_{k(X)}/\Ker (G_{\overline k(X)}\twoheadrightarrow G_{\overline k (X)}^{\ab})$ is the {\it geometrically abelian} 
absolute Galois group of $X$. For $U\subseteq X$ as above we have exact sequences 
$$1\to I_U\to \pi_1(U)^{(\ab)}\to \pi_1(X)^{(\ab)}\to 1,\tag 3$$
where $I_U\defeq \Ker (\pi_1(U)^{(\ab)}\twoheadrightarrow  \pi_1(X)^{(\ab)})=\Ker (\pi_1(U_{\overline k})^{\ab}\twoheadrightarrow 
 \pi_1(X_{\overline k})^{\ab})$, and 
$$1\to I\to G_{k(X)}^{(\ab)}\to \pi_1(X)^{(\ab)}\to 1,\tag 4$$ 
where $I\defeq \Ker (G_{k(X)}^{(\ab)}\twoheadrightarrow  \pi_1(X)^{(\ab)})=\Ker (G_{\overline k(X)}^{\ab}\twoheadrightarrow 
 \pi_1(X_{\overline k})^{\ab})$. Note that $G_{k(X)}^{(\ab)}=\underset {U}\to \varprojlim\ \pi_1(U)^{(\ab)}$, and
 $I=\underset {U}\to \varprojlim\ I_U$, where the limits are over all open subschemes $U\subseteq X$.
Moreover, if $P_1,\ldots,P_n\in X$ are closed points and $U\defeq X\setminus \{P_1,\ldots,P_n\}$ then we have an exact sequence
$$0\to \hat \Bbb Z(1)@>>> \prod_{i=1}^n \Ind_{k(P_i)}^k \hat \Bbb Z(1)\to I_U\to 0,\tag 5$$ 
as follows from the well-known structure of $\pi_1(U_{\overline k},\eta)^{\ab}$, 
and (by passing to the projective limit we obtain) the exact sequence 
$$0\to \hat \Bbb Z(1)\to  \prod_{P\in X^{\cl}} \Ind_{k(P)}^k \hat \Bbb Z(1) \to I\to 0\tag 6$$
of $G_k$-modules, 
where in (6) the product is over all closed points $P\in X^{\cl}$.     
More precisely, for $U=X\setminus \{P_1,\ldots,P_n\}$ as above let $J_U$ be the generalised jacobian of $U$ which sits in the following exact sequence
$$0\to H_U\to J_U\to J\to 0\tag *$$ 
where $H_U\defeq \Coker \left(\Bbb G_{m,k}@>>>\prod_{i=1}^n\Res _{k(P_i)/k}\Bbb G_m\right)$ is a torus and $J\defeq \Jac (X)$ is the jacobian of $X$. 
We have an exact sequence of Tate modules
$$0\to TH_U=I_U\to TJ_U\to TJ\to 0\tag **$$ 
and $TJ_U$ is identified with $\pi_1(U_{\overline k},\eta)^{\ab}$ (as $G_k$-modules).

As was observed in [Esnault-Wittenberg] Remark 2.3(ii), 
in the case where $k$ is a $p$-adic local field, $\index (X)=1$ (i.e., $X$ possesses a divisor of degree $1$) if and only if the exact sequence
(2) splits. Our first main result is the following. 
(See [Lichtenbaum] for the definition of the period of a curve.)

\proclaim {Theorem A} Assume that $k$ is a $p$-adic local field for some prime integer $p\ge 2$ (i.e., $k/\Bbb Q_p$ is a finite extension).
Then there exists an effectively computable non-empty open affine
subscheme $U\subset X$ with the following properties.

{\bf (i)}\ {\bf If} the exact sequence (1) of $\pi_1(U,\eta)^{(\ab)}$ splits then $\period (X)=1$ (i.e., $X$ possesses a $k$-rational divisor class of degree $1$) and $\index (X)$ equals 
$1$ or $2$.

{\bf (ii)} Assume $\period(X)=\index(X)$. Then $\index (X)=1$ (i.e., $X$ possesses a degree 1 divisor) 
{\bf if and only if} the exact sequence (1) of $\pi_1(U,\eta)^{(\ab)}$ splits.
\endproclaim

The term {\it effectively computable} in Theorem A means that one can effectively compute $U$ if one can effectively compute a set of topological generators of 
the group of $k$-rational points $J(k)$ of the jacobian $J\defeq \Jac (X)$ (cf. proof of Theorem A and Lemma 1.1).
For a $p$-adic local field $\ell$, write $(\ell^{\times})^{\wedge}$ for the profinite completion of 
its multiplicative group $\ell^{\times}\defeq \ell\setminus \{0\}$.
Our second main result is the following, in which we compute the torsor of splittings of the exact sequence (2).

\proclaim {Theorem B} With the assumptions in Theorem A, assume that $\index (X)=1$. Then there exists an exact sequence
$$0\to H^1(G_k,I)\to H^1(G_k,G_{\overline k(X)}^{\ab})\to J(k)\to 0,\tag 7$$
as well as isomorphisms $\underset {U}\to \varprojlim \ J_U(k)^{\wedge}\isom H^1(G_k,G_{\overline k(X)}^{\ab})$, $\underset {U}\to \varprojlim \ H_U(k)^{\wedge}\isom H^1(G_k,I)$,
where the projective limit is over all open subschemes $U\subseteq X$, and $J_U(k)^{\wedge}\defeq J_U(k)\otimes_{\Bbb Z} \hat \Bbb Z$.
Moreover, if $U=X\setminus \{P_1,\ldots,P_n\}$ is affine, then we have an exact sequence
$$1\to (k^{\times})^{\wedge}@>>> \prod_{i=1}^n (k(P_i)^{\times})^{\wedge}\to  H^1(G_k,I_U)\to 0,$$
and (by passing to the projective limit we obtain) an exact sequence
$$1\to (k^{\times})^{\wedge}@>>> \prod_{P\in X^{\cl}} (k(P)^{\times})^{\wedge} \to  H^1(G_k,I)\to 0,$$ 
where the product is over all closed points $P\in X^{\cl}$.  
\endproclaim

Next, let $s:G_k\to \pi_1(X,\eta)$ be a {\it section} of the projection $\pi_1(X,\eta)\twoheadrightarrow G_k$. 
Recall that the section $s$ is called {\it orthogonal} to $\Pic^{\wedge}$ (resp. $\Pic^0$) if the homomorphism $s^{\star}:H_{\et}^2(X,\hat \Bbb Z(1))\to H^2(G_k,\hat \Bbb Z(1))$
induced by $s$ [$H^2_{\et}(X,\hat \Bbb Z(1))$ is naturally identified with $H^2(\pi_1(X,\eta),\hat \Bbb Z(1))$ (cf. [Mochizuki], Proposition 1.1)]
annihilates the Picard part $\Pic(X)^{\wedge}\defeq \Pic(X)\otimes _{\Bbb Z}\hat \Bbb Z$ (resp. the (image in $\Pic(X)^{\wedge}$ of the) degree $0$ part $\Pic^0(X)$) 
of $H_{\et}^2(X,\hat \Bbb Z(1))$ (cf. [Sa\"\i di], Definition 1.4.1).
We say that the section $s$ is {\it strongly orthogonal} to $\Pic^{\wedge}$ (resp. $\Pic^0$) if for every neighbourhood $X_i\defeq X_i[s]\to X$ of the section $s$ and the induced section
$s_i:G_k\to \pi_1(X_i,\eta)$ of the projection
$\pi_1(X_i,\eta) \twoheadrightarrow G_k$ (cf. loc. cit. 1.3) the section $s_i$ is orthogonal to $\Pic^{\wedge}$ (resp. $\Pic^0$), $i\ge 1$.
(Note that the above definition differs slightly from the definition in loc. cit. where the notion of having a cycle class orthogonal to $\Pic^{\wedge}$
was defined as being strongly orthogonal to $\Pic^{\wedge}$ in the above sense.) We say that the section $s$ is {\it uniformly orthogonal} to $\Pic^{\wedge}$ (resp. $\Pic^0$) if given a finite extension $\ell/k$ and the induced section $s_{\ell}:G_{\ell}\to \pi_1(X_{\ell},\eta)$ of the projection $\pi_1(X_{\ell},\eta)\twoheadrightarrow G_{\ell}$, where $X_{\ell}\defeq X\times_k\ell$, then $s_{\ell}$ is orthogonal to $\Pic^{\wedge}$ (resp. $\Pic^0$). 
The above definitions carry out in a similar way in the case of sections of geometrically pro-$\Sigma$ arithmetic fundamental groups, where $\Sigma$ 
is a non-empty set of prime integers (cf. loc. cit.).

To a section $s:G_k\to \pi_1(X,\eta)$ as above one associates naturally, by considering the composite morphism of $s$ and the natural projection 
$\pi_1(X,\eta)\twoheadrightarrow \pi_1(X,\eta)^{(\ab)}$, a section $s^{\ab}:G_k\to \pi_1(X,\eta)^{(\ab)}$
of the projection $\pi_1(X,\eta)^{(\ab)}\twoheadrightarrow G_k$.
Let $J^1\defeq \Pic^1_X$ which is a torsor under $J$. There is a natural morphism $X\to J^1$.
In case $\period (X)=1$, hence $J^1(k)\neq \emptyset$,
we identify $J^1$ and $J$ via the isomorphism $J^1\isom J$ which maps a point $z\in J^1(k)$ to the zero section $0\in J(k)$
and consider the composite morphism $X\to J^1\isom J$. We then obtain a commutative diagram
$$
\CD
1@>>> \pi_1(X_{\overline k},\eta)^{\ab}@>>> \pi_1(X,\eta)^{(\ab)}@>>> G_k @>>> 1\\
@. @VVV   @VVV  @VVV  @.\\
1@>>> \pi_1(J_{\overline k},\eta)@>>> \pi_1(J,\eta)@>>> G_k@>>> 1\\
\endCD
$$
where the vertical maps are isomorphisms. We fix compatible base points of the torsors of splittings of the horizontal sequences in the above diagram. 
For example, the splitting $s_z:G_k\to \pi_1(X,\eta)^{(\ab)}$ of the upper sequence arising from the above point $z\in J^1(k)$ 
(once we identify $\pi_1(X,\eta)^{(\ab)}$ and $\pi_1(J^1,\eta)$),
and the induced splitting $s_0:G_k\to \pi_1(J,\eta)$
of the lower sequence which arises from the zero section $0\in J(k)$.
The section $s^{\ab}:G_k\to \pi_1(X,\eta)^{(\ab)}$ gives rise to a section $s^{\ab}:G_k\to \pi_1(J,\eta)$ of the lower sequence in the above diagram, we will denote by 
$[s^{\ab}]\defeq [s^{\ab}-s_0]\in H^1(G_k,TJ)$ the cohomology class (i.e., the cohomology class of the $1$-cocycle $s^{\ab}-s_0:G_k\to \pi_1(J_{\overline k},\eta)$)
associated to $s^{\ab}$, where $TJ$ is the Tate module of $J$  which we identify with $\pi_1(J_{\overline k},\eta)$.
Recall the Kummer exact sequence $0\to J(k)^{\wedge} \to H^1(G_k,TJ)\to TH^1(G_k,J)\to 0$, we view $J(k)^{\wedge}\defeq J(k)\otimes _{\Bbb Z} \hat \Bbb Z$ as a subgroup of $H^1(G_k,TJ)$ via the Kummer map
$J(k)^{\wedge}\to H^1(G_k,TJ)$ (cf. [Sa\"\i di1], $\S1$, for a detailed discussion). 
If $k$ is a $p$-adic local field then the natural map $J(k)\to J(k)^{\wedge}$ is an isomorphism as follows from the well-known structure of $J(k)$ in this case.
In this paper, if $k/\Bbb Q_p$ is a finite extension, we will identify $J(k)$ and $J(k)^{\wedge}$ via this isomorphism.
Our next main result is the following which characterises sections of arithmetic fundamental groups of curves over 
$p$-adic local fields which are orthogonal to $\Pic^0$.

\proclaim {Theorem C} With the assumptions in Theorem A, let $s:G_k\to \pi_1(X,\eta)$ be a section of the projection
$\pi_1(X,\eta) \twoheadrightarrow G_k$. Then the followings hold.

{\bf (i)}\ The section $s$ is orthogonal to $\Pic^0$ {\bf if} (resp. assuming $\index (X)=1$, {\bf if and only if}) the section $s^{\ab}:G_k\to \pi_1(X,\eta)^{(\ab)}$ 
lifts to a section $\tilde s^{\ab}:G_k\to G_{k(X)}^{(\ab)}$ of the exact sequence (2).

{\bf (ii)} Assume that $X(k)\neq \emptyset$. Then $s$ is orthogonal to $\Pic^0$ {\bf if and only if} $[s^{\ab}]\in J(k)$. 
\endproclaim

The assumption that $X(k)\neq \emptyset $ in Theorem C(ii) is rather mild. Indeed, in order to verify that $s$ is orthogonal to $\Pic^{\wedge}$ (resp. $\Pic^0$) one can pass to a finite extension $\ell/k$,
and the corresponding section $s_{\ell}:G_{\ell}\to \pi_1(X_{\ell},\eta)$ of the projection $\pi_1(X_{\ell},\eta)\twoheadrightarrow G_{\ell}$ (cf. proof of Theorem C(i)).
Thus, Theorem C (especially Theorem C(ii)) can be in principle used to detect if a section $s$ as above is (strongly) orthogonal to $\Pic^0$.
As an illustration of this fact we observe that the {\it non-geometric} (geometrically pro-$p$) section constructed by Hoshi over $p$-adic local fields in 
[Hoshi] is orthogonal to $\Pic^0$ (cf. Proposition 3.3). Finally, we observe the following characterisation of sections $s$ as above which are strongly orthogonal to $\Pic^{\wedge}$.

\proclaim {Theorem D} With the assumptions in Theorem A, let $s:G_k\to \pi_1(X,\eta)$ be a section of the projection
$\pi_1(X,\eta) \twoheadrightarrow G_k$. Then the following two conditions are equivalent.

{\bf (i)} The section $s$ is strongly orthogonal to $\Pic^{\wedge}$.

{\bf (ii)} For every neighbourhood $X_i\defeq X_i[s]\to X$ of $s$, $i\ge 1$ (cf. above discussion),
the section $s_i^{\ab}:G_k\to \pi_1(X_i,\eta)^{(\ab)}$ lifts to a section $\tilde s_i:G_k\to G_{k(X_i)}^{(\ab)}$ of the projection $G_{k(X_i)}^{(\ab)}\twoheadrightarrow G_k$.
\endproclaim

\subhead {Acknowledgment} 
\endsubhead I would like to thank the referee for his/her careful reading of the paper and useful comments. I thank Akio Tamagawa for the 
several interesting discussions we had on the topic of this paper, especially around Theorem A.

\subhead {\S 1. Proof of Theorem A}
\endsubhead In this section we prove Theorem A. First, note that if $\index (X)=1$ (i.e., $X$ possesses a divisor of degree $1$) 
then the exact sequence (1) (as well as the exact sequence (2)) splits for every open subscheme $U\subseteq X$,
as follows from a restriction and corestriction argument in Galois cohomology. 
We start with the following Lemmas.

\proclaim {Lemma 1.1} There exists an effectively computable open affine
subscheme $U\subset X$ such that $H^1(G_k,J_U)$ is {\bf finite}.
\endproclaim


\proclaim {Lemma 1.2} There exists an effectively computable open affine
subscheme $U\subset X$ such that $\Ker ( H^1(G_k,J)\to H^2(G_k,H_U))$ is {\bf finite}, where the map $H^1(G_k,J)\to H^2(G_k,H_U)$ arises from (the Galois cohomology of) 
the exact sequence (*).
\endproclaim
\demo{Proof of Lemma 1.2}
Let $U=X\setminus \{P_1,\ldots,P_n\}$ be an open affine subscheme ($P_1,\ldots,P_n\in X$ are closed points).
We have an exact sequence (where $\Br$ denotes Brauer groups) 
$H^1(G_k,H_U)\to \Br(k)\to \oplus _{i=1}^n \Br (k(P_i))\to H^2(G_k,H_U)\to 0$ arising from the long Galois cohomology exact sequence associated to 
the exact sequence  $1\to \Bbb G_{m,k} \to \prod_{i=1}^n\Res _{k(P_i)/k}\Bbb G_m \to H_U\to 1$ of $G_k$-modules [note that by Shapiro's Lemma we have $H^2(G_k,\Res_{k(P_i)/k}\Bbb G_m)\isom H^2(G_{k(P_i)},\Bbb G_m)=\Br(k(P_i)$], and we identify the Brauer group of a $p$-adic local field with 
$\Bbb Q/\Bbb Z$. The Pontryagin dual of
$H^2(G_k,H_U)$ is identified  with $\Div^0(X\setminus U)^{\wedge}\defeq\Ker (\oplus _{i=1}^n\hat \Bbb Z.P_i @>\deg>> \hat \Bbb Z)$ where $\deg(P_i)=[k(P_i):k]$. 
The dual of $H^1(G_k,J)$ is (by Tate duality) $J(k)$, and the dual of the map $H^1(G_k,J)\to H^2(G_k,H_U)$ is the homomorphism  
$\Div^0(X\setminus U)^{\wedge}\to J(k)^{\wedge}$ which is induced by the map $\Div^0(X\setminus U)\to J(k)$ which maps a divisor of degree $0$ on $X$ supported on 
$X\setminus U$ to its class in $J(k)$. 
Further, $J(k)$ is topologically finitely generated as is well-known (cf. [Mattuck]). Let $\{x_1,\ldots,x_t\}$ be topological generators of $J(k)$. There exists an integer $r\ge 1$ 
depending only on $g$ (for example $2$ if $g=1$, or $g-1$ if $g>1$) 
such that  $rx_i=[D_i]$ is the class of a degree $0$ divisor $D_i=\sum_{j=1}^{m_i}n_{i,j}P_{i,j}$ on $X$, for $1\le i\le t$. Now let $U\defeq X\setminus \{P_{i,j}\}_{i=1,j=1}^{t,m_i}$. Then 
$\Im (\Div^0(X\setminus U)^{\wedge}\to J(k))$ has finite index in $J(k)$, and by duality $\Ker (H^1(G_k,J)\to H^2(G_k,H_U))$ is finite. 
\qed
\enddemo

\demo{Proof of Lemma 1.1}
Let $U\subset X$ be as in Lemma 1.2. We have an exact sequence $H^1(G_k,H_U)\to H^1(G_k,J_U)\to H^1(G_k,J)\to H^2(G_k,H_U)$
(arising from the long Galois cohomology exact sequence associated to (*), cf. diagram below). Further, $H^1(G_k,H_U)$ is finite (cf. [Serre], II.5.8 Theorem 6), $\Ker \left (H^1(G_k,J)\to H^2(G_k,H_U)\right )$ is finite (cf. Lemma 1.2), hence $H^1(G_k,J_U)$ is finite as follows from the exactness of the above sequence. This finishes the proof of Lemma 1.1.
\qed
\enddemo

Next, we resume the proof of Theorem A. Let $U\subset X$ be an open affine subscheme. We have a commutative diagram 

$$
\CD
J(k)^{\wedge} @>>> H^1(G_k,TJ)\\
@VVV   @VVV\\
H^1(G_k,H_U) @>>> H^2(G_k,TH_U) \\
@VVV   @VVV\\
H^1(G_k,J_U) @>>> H^2(G_k,TJ_U)\\
@VVV  @VVV\\
H^1(G_k,J) @>>> H^2(G_k,TJ)\\
@VVV   @VVV\\
H^2(G_k,H_U)@>>> H^3(G_k,TH_U)=0\\
\endCD
$$
where the vertical sequences are exact and arise from the exact sequences (*) and (**),
and the horizontal maps are Kummer homomorphisms arising from the Kummer exact sequences in Galois cohomology associated to the algebraic groups $J$, $H_U$, and $J_U$, respectively.
The middle (resp. fourth from the top) horizontal map maps the class $[J_U^1]$ of the universal torsor $J_U^1$ (of degree $1$) (resp. the class $[J^1]$ of $J^1=\Pic^1_{X/k}$)
to the class $[\pi_1(U,\eta)^{(\ab)}]$ of the group extension $\pi_1(U,\eta)^{(\ab)}$ (resp. the class $[\pi_1(X,\eta)^{(\ab)}]$ of the group extension $\pi_1(X,\eta)^{(\ab)}$)
(this is a well-known fact, see for example [Harari-Szamuely] Proposition 2.2 and Remark 2.4).
Further, $[J_U^1]$ (resp. $[\pi_1(U,\eta)^{(\ab)}]$) maps to $[J^1]$ (resp. $[\pi_1(X,\eta)^{(\ab)}]$) under the left third vertical map from the top 
(resp. right third vertical map from the top).

Next, we let $U$ be as in Lemma 1.1. We prove that assertions (i) and (ii) in Theorem A are satisfied in this case.

We prove assertion (i). Assume that the class $[\pi_1(U,\eta)^{(\ab)}]$ is trivial in $H^2(G_k,TJ_U)$ which implies that the class $[J_U^1]$ is divisible in
$H^1(G_k,J_U)$. 
(The map $H^1(G_k,J_U)\to H^2(G_k,TJ_U)$ factors through $\underset {n} \to \varprojlim\ H^1(G_k,J_U)/nH^1(G_k,J_U$) and the latter group injects into  $H^2(G_k,TJ_U)$). 
As the group $H^1(G_k,J_U)$ is finite the class of $[J^1_U]$ is then trivial. Thus, $[J^1]=0$ in  $H^1(G_k,J)$ (cf. above discussion) which implies that $X$ possesses a $k$-rational divisor class of degree $1$, i.e., $\period(X)=1$. The rest of the assertion follows from the fact that either $\index(X)=\period (X)$ or $\index (X)=2\period (X)$ (cf. [Lichtenbaum], Theorem 7).

Assertion (ii) follows from (i) for the if part, and the only if part follows from the observation at the start of the proof of Theorem A. This finishes the proof of Theorem A.
\qed

\subhead {\S 2. Proof of Theorem B}
\endsubhead
In this section we prove Theorem B. We use the same assumptions as in Theorem A and further suppose that $X$ possesses a degree $1$ divisor. 
We start with the following lemma.

\proclaim{Lemma 2.1} We use the assumptions in Theorem A. Assume that $\index (X)=1$. Then $\underset {U}\to \varprojlim \ H^1(G_k,H_U)=0$
where the limit is over all non-empty open subschemes of $X$.

\endproclaim

\demo{Proof of Lemma 2.1} The exact sequence $1\to \Bbb G_{m,k}\to \prod_{i=1}^n\Res _{k(P_i)/k}\Bbb G_m\to H_U\to 1$ induces in cohomology an exact sequence
$0\to  H^1(G_k,H_U)\to \Br(k)\to \prod_{i=1}^n \Br (k(P_i))$ (note that $H^1(G_k,\Res _{k(P_i)/k}\Bbb G_m)\isom H^1(G_{k(P_i)},\Bbb G_m)=0$, 
and $H^2(G_k,\Res _{k(P_i)/k}\Bbb G_m)\isom H^2(G_{k(P_i)},\Bbb G_m)=\Br(k(P_i))$, 
by Shapiro's Lemma)
and by passing to the projective limit over all $U\subset X$ open we obtain an exact sequence
$0\to \underset {U}\to \varprojlim \ H^1(G_k,H_U)\to \Br(k)\to \prod_{P\in X^{\cl}} \Br (k(P))$. Now $\Ker (\Br(k)\to \prod_{P\in X^{\cl}} \Br (k(P)))$
is finite of cardinality $\index (X)$ (cf. [Lichtenbaum], Theorem 3), which equals $1$ under our assumption that $X$ possesses a degree one divisor.
\qed
\enddemo

Next, we resume the proof of Theorem B. Consider the morphism $X\to J$ as in the introduction, and identify the $G_k$-modules $TJ$ and $\pi_1(X_{\overline k})^{\ab}$.
The assertions regarding the structure of $H^1(G_k,I_U)$ and $H^1(G_k,I)$ follow 
easily from Kummer theory (consider the long cohomology exact sequences associated to the exact sequences (5) and (6) of $G_k$-modules). 
We establish the exact sequence (7) in the statement of the theorem as well as the isomorphisms  $\underset {U}\to \varprojlim \ J_U(k)^{\wedge}\isom H^1(G_k,G_{\overline k(X)}^{\ab})$
and  $\underset {U}\to \varprojlim \ H_U(k)^{\wedge}\isom H^1(G_k,I)$ therein.

We have a commutative diagram of group homomorphisms
$$
\CD   
@.   @.   0   @.  0    @.  0\\
@.    @.   @AAA    @AAA   @AAA   @.\\
@. TH^1(G_k,H_U)=0  @>>> TH^1(G_k,J_U) @>>> TH^1(G_k,J) @>>> TH^2(G_k,H_U)\\
@.   @AAA  @AAA    @AAA   @AAA\\
0 @>>> H^1(G_k,TH_U)  @>>> H^1(G_k,TJ_U) @>>> H^1(G_k,TJ) @>>> H^2(G_k,TH_U)\\
@.   @AAA  @AAA @AAA   @AAA\\
0 @>>> H_U(k)^{\wedge}  @>>> J_U(k)^{\wedge}  @>>> J(k)^{\wedge} @>>> H^1(G_k,H_U)^{\wedge}\\
@.   @AAA  @AAA @AAA   @AAA\\
@. 0  @.  0  @. 0 @. 0\\
\endCD
$$
where the vertical sequences are Kummer exact sequences, and the middle and lower horizontal sequences arise from the exact sequences (*) and (**). 
Note that since $H^1(G_k,H_U)$ is finite (cf. [Serre], II.5.8 Theorem 6), $TH^1(G_k,H_U)=0$, and the natural map $H^1(G_k,H_U)\to H^1(G_k,H_U)^{\wedge}$ is an  isomorphism. 
The middle horizontal sequence is exact and arises from the long cohomology exact sequence associated to the exact sequence (**). (Note that $H^0(G_k,TJ)=0$ as follows from the well-known fact that $J(k)^{\tor}$ is finite.) 

The map $TH^1(G_k,J_U) \to TH^1(G_k,J)$ is injective as follows easily from the exact sequence $H^1(G_k,H_U)\to H^1(G_k,J_U) \to H^1(G_k,J)$, the left exactness of the inverse limit functor, and the fact that $H^1(G_k,H_U)$ is finite (cf. [Serre], II.5.8 Theorem 6). We claim that the lower horizontal sequence is exact. Indeed, the map $H_U(k)^{\wedge}  \to J_U(k)^{\wedge}$ is injective as follows from the commutativity of the far left lower square, and the injectivity of the maps   $H_U(k)^{\wedge}\to H^1(G_k,TH_U) \to H^1(G_k,TJ_U)$. Exactness at $J_U(k)^{\wedge}$ follows from the commutativity of the lower middle square,  the exactness at $H^1(G_k,TJ_U)$ of the middle horizontal exact sequence, and the fact that the map $H_U(k)^{\wedge}\to H^1(G_k,TH_U)$ is an isomorphism.
Let $\alpha\in J(k)^{\wedge}$ with trivial image in $H^1(G_k,H_U)^{\wedge}$, its image $\alpha \in H^1(G_k,TJ)$ is the image of an element $\beta \in  H^1(G_k,TJ_U)$ by the commutativity of the right lower square and the exactness of the middle horizontal sequence. As $\alpha$ maps to $0$ in $TH^1(G_k,J)$, the image of $\beta$ in $TH^1(G_k,J_U)$ is $0$ by the commutativity of the middle upper square and the injectivity of the map $TH^1(G_k,J_U) \to TH^1(G_k,J)$. Thus $\beta \in  J_U(k)^{\wedge}$ maps to $\alpha$ in  $J(k)^{\wedge}$ and the lower sequence is 
exact at $J(k)^{\wedge}$.  

By passing to the projective limit over all open subschemes $U\subseteq X$ we obtain a commutative diagram 
$$
\CD
@.     @.    @.  0  @. 0\\
@.    @.     @. @AAA   @AAA\\
@.   0  @>>> \underset {U}\to \varprojlim\ TH^1(G_k,J_U) @>>> TH^1(G_k,J) @>>> \underset {U}\to \varprojlim\ TH^2(G_k,H_U)\\
@.     @AAA    @AAA   @AAA    @AAA\\
0 @>>> H^1(G_k,I)  @>>> H^1(G_k,G_{\overline k(X)}^{\ab}) @>>> H^1(G_k,TJ) @>>> \underset {U}\to \varprojlim\ H^2(G_k,TH_U)\\
@.   @AAA  @AAA @AAA    @AAA \\
0 @>>> \underset {U}\to \varprojlim \ H_U(k)^{\wedge}  @>>> \underset {U}\to \varprojlim \ J_U(k)^{\wedge}  @>>> J(k) @>>> \underset {U}\to \varprojlim \ H^1(G_k,H_U) \\
@.   @AAA  @AAA @AAA   @AAA\\
@. 0  @.  0  @. 0 @. 0\\
\endCD
$$
where the middle horizontal sequence is exact and arises from the long exact cohomology sequence associated to the exact sequence (4). The left vertical map is an isomorphism
[$H^1(G_k,I)\isom \underset {U}\to \varprojlim\ H^1(G_k,TH_U)$], the second left vertical sequence is exact as follows from the left exactness of the inverse limit functor,
the second right vertical sequence is the Kummer exact sequence associated to $J$, and the right vertical sequence is exact since the $H^1(G_k,H_U)$ are finite (cf. loc. cit.),
thus the Mittag-Leffler condition is satisfied. The map $\underset {U}\to \varprojlim\ TH^1(G_k,J_U) \to TH^1(G_k,J)$ is injective, and
the lower horizontal sequences is exact as follows easily from the left exactness of the inverse limit functor and a similar argument as the one used for the previous diagram 
for the exactness at $J(k)$.

Now,  $\underset {U}\to \varprojlim\ TH^1(G_k,J_U)=0$, as
$\underset {U}\to \varprojlim\ TH^1(G_k,J_U)$ is identified with the intersection of the images of 
$TH^1(G_k,J_U)$ in $TH^1(G_k,J)$, and Lemma 1.1 implies the existence of $U\subset X$ open affine such that $H^1(G_k,J_U)$ is finite hence $TH^1(G_k,J_U)=0$. 
This implies that the injective map $\underset {U}\to \varprojlim \ J_U(k)^{\wedge} \to H^1(G_k,G_{\overline k(X)}^{\ab})$
is an isomorphism and $\Im (H^1(G_k,G_{\overline k(X)}^{\ab})\to H^1(G_k,TJ))$ 
is contained in $J(k)$ (we identify the latter with its image via the injective Kummer map $J(k)\hookrightarrow H^1(G_k,TJ)$).
Further, $\underset {U}\to \varprojlim \ H^1(G_k,H_U)=0$ by Lemma 2.1.
Hence $\Im (H^1(G_k,G_{\overline k(X)}^{\ab})\to H^1(G_k,TJ))=J(k)$ and we obtain the exact sequence (7) as claimed in Theorem B.
This finishes the proof of Theorem B.
\qed

\definition{Remark 2.2} Let $\Sigma$ be a non-empty set of prime integers. The same proof as above yields a pro-$\Sigma$ analog of Theorem B. More precisely, 
let $G_{\overline k(X)}^{\ab,\Sigma}$ (resp. $\pi_1(X_{\overline k},\eta)^{\ab,\Sigma}$)
be the maximal pro-$\Sigma$ quotient of  $G_{\overline k(X)}^{\ab}$ (resp. $\pi_1(X_{\overline k},\eta)^{\ab}$)
which sits in the exact sequence
$0\to I_{\Sigma}\to G_{\overline k(X)}^{\ab,\Sigma} \to \pi_1(X_{\overline k},\eta)^{\ab,\Sigma} \to 0$,
where $I_{\Sigma}\defeq \Ker (G_{\overline k(X)}^{\ab,\Sigma} \twoheadrightarrow  \pi_1(X_{\overline k},\eta)^{\ab,\Sigma})$.
Then, with the same assumptions as in Theorem B, we have an exact sequence
$0\to H^1(G_k,I_{\Sigma})\to H^1(G_k,G_{\overline k(X)}^{\ab,\Sigma})\to J(k)^{\Sigma}\to 0$,
where $J(k)^{\Sigma}$ is the maximal pro-$\Sigma$ quotient of $J(k)$.
\enddefinition

\subhead {\S 3. Proof of Theorem C} 
\endsubhead
In this section we prove Theorem C, we use the same assumptions as in Theorem A. The following Lemma will be useful.

\proclaim{Lemma 3.1} Let $s:G_k\to \pi_1(X,\eta)$ be a {\it section} of the projection $\pi_1(X,\eta)\twoheadrightarrow G_k$. If $s$ is orthogonal to $\Pic^0$
then $s$ is uniformly orthogonal to $\Pic^0$.
\endproclaim

\demo{Proof} Similar to the proof of Proposition 1.6.7 in [Sa\"\i di].
\qed
\enddemo

\demo {Proof of Theorem C(i)} 
First, assume that $s^{\ab}:G_k\to \pi_1(X,\eta)^{(\ab)}$ lifts to a section $\tilde s^{\ab}:G_k\to G_{k(X)}^{(\ab)}$ of the exact sequence (2).
We will show that $s$ is orthogonal to $\Pic^0$. Let $\Cal L\in \Pic^0(X)$ corresponding to the class of a degree zero divisor $D=\sum_{i=1}^tn_iP_i$. 
Given a finite extension $\ell/k$, $X_{\ell}\defeq X\times _k\ell$, we have a commutative diagram
$$
\CD
\Pic^0(X_{\ell}) @>>> H^2(X_{\ell},\hat \Bbb Z(1)) @>{s_{\ell}^{\star}}>> H^2(G_{\ell},\hat \Bbb Z(1))\\
@AAA   @AAA  @AAA\\
\Pic^0(X) @>>> H^2(X,\hat \Bbb Z(1)) @>{s^{\star}}>> H^2(G_k,\hat \Bbb Z(1))\\
\endCD
$$ 
where the left lower and upper horizontal maps arise from Kummer theory (they are injective), the vertical maps are restriction maps, and the map $s_{\ell}^{\star}$ is induced by the section $s_{\ell}:G_k\to 
\pi_1(X_{\ell},\eta)$ of the projection $\pi_1(X_{\ell},\eta)\twoheadrightarrow G_k$ which is induced by $s$. 
Identifying both $H^2(G_k,\hat \Bbb Z(1))$ and $H^2(G_{\ell},\hat \Bbb Z(1))$
with $\hat \Bbb Z$, the far right vertical map is multiplication by the degree $[\ell:k]$ of $\ell/k$. In particular, this map is injective. To show that the image of $\Cal L$ in 
$H^2(G_k,\hat \Bbb Z(1))$ is trivial it thus suffices to show that its image in $H^2(G_{\ell},\hat \Bbb Z(1))$ is trivial. We can then, without loss of generality,
and after possibly pulling back the line bundle $\Cal L$ to $X_{\ell}$ for a suitable finite extension $\ell/k$, assume that the points $P_1,\ldots,P_t\in X$ are $k$-rational and $\deg(\Cal L)=\sum_{i=1}^tn_i=0$. Let $U\defeq X\setminus \{P_1,\ldots,P_t\}$.

Consider the following commutative diagram of horizontal exact sequences.
$$
\CD
@.    \hat \Bbb Z(1)    @= \hat \Bbb Z(1)\\
@.   @V{\diag}VV   @V{\diag}VV\\
1 @>>> I_U^{\cn} @>>> \pi_1(U,\eta)^{\c-\cn} @>>>  \pi_1(X,\eta) @>>> 1\\
@. @VVV @VVV  @V{\id}VV\\
1 @>>> I_U @>>> \widetilde \pi_1(U,\eta) @>>>  \pi_1(X,\eta) @>>> 1\\
@. @V{\id}VV @VVV  @VVV\\
1@>>> I_U@>>> \pi_1(U,\eta)^{(\ab)}@>>> \pi_1(X,\eta)^{(\ab)}@>>> 1\\
@. @VVV  @VVV @VVV\\
@.  1   @. 1 @. 1\\
\endCD
$$
Here the group extension $\widetilde \pi_1(U,\eta)$ is the pull back of the lower horizontal exact sequence by the map $\pi_1(X,\eta)\twoheadrightarrow \pi_1(X,\eta)^{(\ab)}$
(i.e., the lower right square is cartesian), 
$\pi_1(U,\eta)^{\c-\cn}$
is the geometrically cuspidally central quotient of $\pi_1(U,\eta)$ (cf. [Sa\"\i di] 2.1.1), the surjective map $\pi_1(U,\eta)^{\c-\cn}\twoheadrightarrow 
\widetilde \pi_1(U,\eta)$ is the natural one 
($\pi_1(X_{\overline k},\eta)$ acts trivially on $I_U$), $I_U^{\cn}$ is the $G_k$-module $\prod_{i=1}^t\hat \Bbb Z(1)$ (cf. loc. cit. proof of Lemma 2.3.1), $\hat \Bbb Z(1)@>\diag>>
I_U^{\cn}=\prod_{i=1}^t\hat \Bbb Z(1)$ is the diagonal embedding, and we have an exact sequence of $G_k$-modules
$0\to \hat \Bbb Z(1)@>{\diag}>> I_U^{\cn}=\prod_{i=1}^t\hat \Bbb Z(1)\to I_U\to 0$.

By pulling back the group extension $\pi_1(U,\eta)^{\c-\cn}$ by the section $s:G_k\to \pi_1(X,\eta)$ we obtain a group extension
$1\to I_U^{\cn}\to F_U\to G_k\to 1$. Further, by pulling back the group extension $1\to I_U\to \pi_1(U,\eta)^{(\ab)}\to \pi_1(X,\eta)^{(\ab)}\to 1$
by the section $s^{\ab}$ we obtain a group extension 
$1\to I_U\to E_U\to G_k\to 1$, which splits since by assumption $s^{\ab}$ lifts to a section $s_U^{\ab}:G_k\to \pi_1(U,\eta)^{(\ab)}$ of the exact sequence (1). (More precisely, the 
section $s_U^{\ab}$ is induced by $\tilde s^{\ab}$.)
Consider the Galois cohomology exact sequence $H^2(G_k,\hat \Bbb Z(1)) @>{\diag}>> H^2(G_k,I_U^{\cn})=\prod _{i=1}^tH^2(G_k,\hat \Bbb Z(1))\to H^2(G_k,I_U)\to 0$.
The class of the extension $F_U$ in $H^2(G_k,I_U^{\cn})$ coincides with $(s^{\star}(\Cal O(P_i))_{i=1}^t$ (cf. [Sa\"\i di] proof of Lemma 2.3.1), and
the class of the group extension $E_U$ in $H^2(G_k,I_U)$ is the image of the class of $F_U$ via the above map 
$H^2(G_k,I_U^{\cn})\to H^2(G_k,I_U)$. In particular, since the class of $E_U$ vanishes in $H^2(G_k,I_U)$,
the class of $F_U$ lies in the diagonal image of $H^2(G_k,\hat \Bbb Z(1))$.
Thus, we deduce that $s^{\star}(\Cal O(P_i))$ is independent of $1\le i\le t$ (i.e., equals the same element of $H^2(G_k,\hat \Bbb Z(1))$), and $s^{\star}(\Cal L)=0$. 

Next, we show that the converse holds assuming $\index (X)=1$. We assume that $s$ is orthogonal to $\Pic^0$, $\index (X)=1$, and show that the section $s^{\ab}:G_k\to \pi_1(X,\eta)^{(\ab)}$ lifts to a section $\tilde s^{\ab}:G_k\to G_{k(X)}^{(\ab)}$ of the exact sequence (2).
Recall the exact sequence $1\to I\to G_{k(X)}^{(\ab)}\to \pi_1(X,\eta)^{(\ab)}\to 1$ (resp. $1\to I_U\to \pi_1(U,\eta)^{(\ab)}\to \pi_1(X,\eta)^{(\ab)}\to 1$, for 
$U\subseteq X$ open). By pulling back this group extension by the section $s^{\ab}$ we obtain a group extension 
$1\to I\to E\to G_k\to 1$ (resp. $1\to I_U\to E_U\to G_k\to 1$, for $U\subseteq X$ open), we will show that the group extension $E$ 
is a split extension which would imply the above assertion. Note that $E=\underset {U}\to \varprojlim \ E_U$.

We have a natural identification $H^2(G_k,I)\isom \underset {U}\to \varprojlim\ H^2(G_k,I_U)$, where the limit is over all $U\subseteq X$ as above. Further, for $U\subseteq X$ as above, we have a Kummer exact sequence
$0\to H^1(G_k,H_U)\to H^2(G_k,I_U)\to TH^2(G_k,H_U)\to 0$ (cf. far right vertical sequence in the first diagram in the proof of Theorem B and the identification $I_U\isom TH_U$ of $G_k$-modules), and 
by passing to the projective limit over all $U$ we obtain an exact sequence $0\to \underset {U}\to \varprojlim\ H^1(G_k,H_U)\to \underset {U}\to \varprojlim\ H^2(G_k,I_U)\to 
\underset {U}\to \varprojlim\ TH^2(G_k,H_U)\to 0$, hence an identification $H^2(G_k,I)\isom \underset {U}\to \varprojlim\ H^2(G_k,I_U)\isom \underset {U}\to 
\varprojlim\ TH^2(G_k,H_U)$
since $\underset {U}\to \varprojlim\ H^1(G_k,H_U)=0$ if $\index (X)=1$ (cf. Lemma 2.1). Write $\Tilde E_U$ for the image of the class of the group extension $E_U$ in $TH^2(G_k,H_U)$ via the above map
$H^2(G_k,I_U)\to TH^2(G_k,H_U)$. We will show $\Tilde E_U=0$, $\forall U\subseteq X$ as above, from which it will follow that the class of the group extension $E$ in $H^2(G_k,I)$ 
is trivial.

Let $U=X\setminus \{P_1,\ldots,P_t\}$ be an open affine subscheme, and $k'/k$ a finite extension which splits the torus $H_U$. We have the following commutative diagram 
of Kummer exact sequences
$$
\CD
0 @>>> H^1(G_{k'},H_U)=0@>>> H^2(G_{k'},I_U) @>>>  TH^2(G_{k'},H_U)@>>> 0\\
@. @A{\res}AA    @A{\res}AA  @A{\res}AA \\
0 @>>> H^1(G_{k},H_U)@>>> H^2(G_{k},I_U) @>>>  TH^2(G_{k},H_U)@>>> 0\\
\endCD
$$
where the vertical maps are restriction maps. We claim $\res ([E_U])=0$ in $H^2(G_{k'},I_U)$.
Indeed, first, using Lemma 3.1, we can (without loss of generality) assume that $\{P_1,\ldots,P_t\}\subset X(k)$, $k=k'$, and we have to show that $[E_U]=0$. 
Recall the commutative diagram and notations in the proof of the if part of Theorem C(i) above, as well as the discussion therein. The assumption that $s$ is orthogonal to $\Pic^0$ implies, by considering the classes of the various degree zero divisors $P_i-P_j$ with $1\le i\le t$, $1\le j\le t$, that 
$s^{\star}(\Cal O(P_i))$ is independent of $1\le i\le t$, which implies that the class of $F_U$ lies in the diagonal image of $H^2(G_k,\hat \Bbb Z(1))$ and the class of 
$E_U$ is trivial (cf. loc. cit.).  Thus, $\res(E_U)=0$ in $H^2(G_{k'},I_U)$ as claimed which implies
$\res (\Tilde E_U)=0$ in $TH^2(G_{k'},H_U)$.
Finally, the far right vertical map in the above diagram is injective, from which it follows that $\Tilde E_U=0$.
Indeed, we have a commutative diagram
$$
\CD
H^2(G_{k'},\Bbb G_m) @>>>  \prod_{i=1}^tH^2(G_{k'},\Ind _{k(P_i)}^k\Bbb G_m) @>>> H^2(G_{k'},H_U)@>>> 0\\
@A{\res}AA    @A{\res}AA  @A{\res}AA \\
H^2(G_{k},\Bbb G_m) @>>>  \prod_{i=1}^tH^2(G_{k},\Ind _{k(P_i)}^k\Bbb G_m) @>>> H^2(G_{k},H_U)@>>> 0 \\
\endCD
$$
where the left vertical map is the map $\Br(k)=\Bbb Q/\Bbb Z\to \Br(k')=\Bbb Q/\Bbb Z$ of multiplication by the degree $[k':k]$ of the extension $k'/k$ which has trivial cokernel
(we identify the Brauer group of a $p$-adic local field with $\Bbb Q/\Bbb Z$),
and the middle vertical map has finite kernel, 
from which it follows that $\Ker (H^2(G_{k},H_U)@>{\res}>>H^2(G_{k'},H_U))$ is finite. By passing to Tate modules we deduce that  
the map $TH^2(G_{k},H_U)@>{\res}>>TH^2(G_{k'},H_U)$ is injective as claimed. 

This finishes the proof of Theorem C(i).
\qed
\enddemo

\demo{Proof of Theorem C(ii)} Assume that $X(k)\neq \emptyset$ and let $x\in X(k)$. Recall the discussion in the introduction after the statement of Theorem B.
In the following argument we use the isomorphism $J^1\isom J$ as in loc. cit. arising from $z\defeq \Cal O(x)\in J^1(k)$. The section 
$s_z:G_k\to \pi_1(X,\eta)^{\ab}\isom \pi_1(J,\eta)$ lifts in this case to a section
$\tilde s_z:G_k\to G_{k(X)}^{(\ab)}$ of the exact sequence (2) since $s_z$ arises from a rational point of $X$. We fix compatible base points of the torsors of splittings of the exact sequences
(2) and (1) with $U=X$ associated to the sections $\tilde s_z$ and $s_z$, respectively. Assertion (ii) follows then from assertion (i) and the exact sequence (7) in Theorem B. 
\qed
\enddemo

\definition {Remark 3.2} Let $\Sigma$ be a non-empty set of prime integers.
Similar proofs as above yield pro-$\Sigma$ analogs of Theorem C(i)(ii) (cf. Remark 2.2).
\enddefinition

In [Hoshi] Hoshi constructed an example of a smooth, geometrically connected, hyperbolic curve $X$ over a $p$-adic local field $k$
and a section $\widetilde s:G_k \to \pi_1(X,\eta)^{(p)}$ of the projection $\pi_1(X,\eta)^{(p)}\twoheadrightarrow G_k$,
where  $\pi_1(X,\eta)^{(p)}$ is the geometrically pro-$p$ quotient of $\pi_1(X,\eta)$, which is not geometric, i.e., the section $\widetilde s$ doesn't arise from a $k$-rational point
of $X$. (There is an exact sequence
$1\to \pi_1(X_{\overline k},\eta)^p\to \pi_1(X,\eta)^{(p)}\to \Gal (\overline k/k)\to 1$, and
$\pi_1(X,\eta)^{(p)}$ is obtained by push-out of the exact sequence
$1\to \pi_1(X_{\overline k},\eta)\to \pi_1(X,\eta)\to \Gal (\overline k/k)\to 1$ by the maximal pro-$p$ quotient
$\pi_1(X_{\overline k},\eta)\twoheadrightarrow \pi_1(X_{\overline k},\eta)^p$ of $\pi_1(X_{\overline k},\eta)$.)
Our next observation is the following.

\proclaim{Proposition 3.3} With the above notations, the section $\widetilde s$ is orthogonal to $\Pic^0$.
\endproclaim

\demo{Proof} Indeed in Hoshi's construction it occurs that $X(k)\neq \emptyset$, and $[(\widetilde s)^{\ab}]\in J(k)^{p}$,
where $J(k)^p$ is the maximal pro-$p$ quotient of $J(k)$ (cf. loc. cit., especially Theorem 3.5 and Corollary 3.6).
Thus, the statement follows from the pro-$p$ analog of Theorem C(ii) (cf. Remarks 2.2 and 3.2).
\qed
\enddemo

\definition {Remark 3.4} The author doesn't know at the time of writing this paper (and is interested to know) if the section $\widetilde s$ in Proposition 3.3,
constructed by Hoshi, 
is strongly orthogonal to $\Pic^0$.
\enddefinition

\subhead
{\S 4. Proof of Theorem D}
\endsubhead
Next, we prove Theorem D. First, assume that assertion (ii) in Theorem D holds, we prove that assertion (i) holds, i.e., $s$ is strongly orthogonal to $\Pic^{\wedge}$.
We have a commutative diagram
$$
\CD
\Pic(X_{i+1})   @>{s_{i+1}^{\star}}>> H^2(G_k,\hat \Bbb Z(1))\\
@AAA   @A{\id}AA\\
\Pic(X_{i})   @>{s_i^{\star}}>> H^2(G_k,\hat \Bbb Z(1))\\
\endCD
$$
where the left vertical map is the pull back of line bundles via the finite morphism $X_{i+1}\to X_i$. Let $M_i\defeq \Im \left(\Pic(X_{i})@>{\deg}>> \Bbb Z\right)$ which is a free 
$\Bbb Z$-module of rank 1 with generator $e_i$, for $i\ge 1$. 
The map $\Pic(X_{i})   @>{s_i^{\star}}>> H^2(G_k,\hat \Bbb Z(1))$ factorises as $\Pic(X_{i})   @>{\deg}>> M_i@>{\rho_i}>> H^2(G_k,\hat \Bbb Z(1))$,
since $s$ is strongly orthogonal to $\Pic^0$ by assumption and Theorem C(i),
and we have a commutative diagram
$$
\CD
M_{i+1}   @>{\rho_{i+1}}>> H^2(G_k,\hat \Bbb Z(1))\\
@AAA   @A{\id}AA\\
M_{i}   @>{\rho_i}>> H^2(G_k,\hat \Bbb Z(1))\\
\endCD
$$
where the left vertical map is defined by $e_i\mapsto [X_{i+1}:X_i]e_{i+1}$ and $[X_{i+1}:X_i]$ is the degree of the finite morphism
$X_{i+1}\to X_i$. Thus $\rho_i(e_i)$ is infinitely divisible in $H^2(G_k,\hat \Bbb Z(1))$, hence $\rho_i(e_i)=0$, for all $i\ge 1$, since $H^2(G_k,\hat \Bbb Z(1))\isom \hat \Bbb Z$. This shows  
that the map $\Pic(X_{i}) @>{s_i^{\star}}>> H^2(G_k,\hat \Bbb Z(1))$ is the zero map for all $i\ge 1$ as required.

Conversely, assume that assumption (i) holds, i.e., the section $s$ is strongly orthogonal to $\Pic^{\wedge}$. Then the section $s$ has a cycle class uniformly orthogonal to 
$\Pic^{\wedge}$ in the sense of [Sa\"\i di], Definition 1.4.1 (cf. loc. cit. Proposition 1.6.7). Assertion (ii) follows then from [Sa\"\i di] Theorem 2.3.5 applied to each $X_i$, $i\ge 1$ (note that with the notations of loc. cit. $G_{k(X_i)}^{(\ab)}$ is a quotient of $G_{k(X_i)}^{\c-\ab}$).

This finishes the proof of Theorem D.
\qed

\definition {Remark 4.1} Let $\Sigma$ be a non-empty set of prime integers.
Similar proofs as above yield a pro-$\Sigma$ analog of Theorem D.
\enddefinition

$$\text{References.}$$


\noindent
[Esnault-Wittenberg] Esnault, H., Wittenberg, O., On abelian birational sections, Journal of the American Mathematical Society, 
Volume 23, Number 3 (2010), 713-724.

\noindent
[Harari-Szamuley], Harari, D., Szamuely, T., Galois sections for abelianized fundamental groups, Math. Ann., 344 (2009), 779-800.

\noindent
[Hoshi] Hoshi, Y., Existence of nongeometric pro-$p$ Galois sections of hyperbolic curves, Publ. Res. Inst. Math. Sci. 46 (2010), no. 4, 829-848.

\noindent
[Lichtenbaum] Lichtenbaum, S., Duality theorems for curves over $p$-adic fields,  Invent. Math.  7  (1969), 120--136.

\noindent
[Mattuck] Mattuck, A., Abelian varieties over P-adic ground fields, Annals of Mathematics, Vol 62, No.1 (1955), 92-119.

\noindent
[Mochizuki] Mochizuki, S., Absolute anabelian cuspidalizations of proper hyperbolic curves,  J. Math. Kyoto
Univ.  47  (2007),  no. 3, 451--539.

\noindent
[Sa\"\i di] Sa\"\i di, M., The cuspidalisation of sections of arithmetic fundamental groups, Advances in Mathematics 230
(2012) 1931-1954.

\noindent
[Sa\"\i di1] Sa\"\i di, M., On the section conjecture over finitely generated fields, Publ. RIMS Kyoto Univ, 52 (2016), 335-357.


\noindent
[Serre] Serre, J.-P., Cohomologie galoisienne, Seconde \'edition, 
Lecture Notes in Mathematics, 5, Springer-Verlag, Berlin-Heidelberg-New York, 1962/1963.

\bigskip

\noindent
Mohamed Sa\"\i di

\noindent
College of Engineering, Mathematics, and Physical Sciences

\noindent
University of Exeter

\noindent
Harrison Building

\noindent
North Park Road

\noindent
EXETER EX4 4QF

\noindent
United Kingdom

\noindent
M.Saidi\@exeter.ac.uk

\end
\enddocument